\documentclass[12pt]{article}
\usepackage{amsfonts,amsthm,amssymb,amsmath}
\newcommand{\om}{\omega}

\newcommand{\ty}{\infty}
\newcommand{\di}{\displaystyle}
\newcommand{\va}{\varphi}

\newcommand{\ga}{\gamma}

\newcommand{\na}{\nabla}

\newcommand{\ld}{\ldots}

\newcommand{\noa}{\noalign{\medskip}}

\newcommand{\su}{\subset}
\newcommand{\qu}{\quad}
\newcommand{\fo}{\forall}
\newcommand{\al}{\alpha}

\newcommand{\ep}{\varepsilon}
\newcommand{\pa}{\partial}
\newcommand{\ti}{\times}

\title{Almost Coquaternion Structure}
\author{Constantin Udri\c ste,\\Department of Mathematics IV,\\ Polytechnic Institute,\\ Bucharest, Romania}

\date{}

\begin{document}
\maketitle

Abstract of the PhD Thesis in Mathematical Sciences, prepared under Prof. Dr. Gheorghe Th. Gheorghiu's 
guidance and defended at the Faculty of Mathematics and Mechanics of the Babe\c s-Bolyai University
of Cluj, Romania, on September 28, 1971.

\bigskip

\centerline{\bf Preface}

\bigskip

For some $(2m+2)$-dimensional differentiable manifolds one can define the notion of 
{\it almost complex structure} and for some $(4n+4)$-dimensional diffe-rentiable manifolds 
one can define the notion of {\it almost quaternion structure} by means of three 
almost complex structures satisfying certain conditions. These notions were 
first introduced by C. Ehresmann [7] in 1947 and then were studied by 
other geometricians [5], [6], [17], [20], [24], [25], [30], [48], [58], [59] etc, 
respectively [2], [22], [24], [26], [27], [28], [29], [57] etc.

For some $(2m+1)$-dimensional differentiable manifolds one can define 
the notion of {\it almost cocomplex structure} which may be considered as 
analogous to the almost complex structure  (for even dimensional manifolds). 
This structure was first introduced in 1960 by S. Sasaki and called a 
$(\phi,\xi,\eta)$-structure [34]. Later on, S. Sasaki and his students used the 
name {\it almost contact structure} for $(\phi,\xi,\eta)$, as the structural group of 
the tangent bundle of a manifold which has a $(\phi,\xi,\eta)$-structure may be 
reduced to $U(n) \times 1$ and conversely [5], [18], [19], [35], [41], [42], [58] etc. 
The same structure was studied in the papers of J. Bouzon [3], K. Ogiue and 
M. Okumura [31], K. Ogiue [30] etc, which named it an almost cocomplex structure.

The odd dimensional sphere $S^{2m+1}$ is a typical example of differentiable 
manifold which admits an almost cocomplex structure $(\phi,\xi,\eta)$ [35]. 
We mention that by means of the 1-form $\eta$ which appears in the 
almost cocomplex structure of $S^{2m+1}$, Gh. Vr\u anceanu [55] and 
C. Teleman [43] defined and studied on $S^{2m+1}$ a nonholonomic manifold 
$V^{2m}_{2m+1}$ whose Riemannian metric is that of a symmetric space $V_{2m}$. 
Taking this idea into account, we observed that on the almost cocomplex manifold a 
nonholonomic (holonomic) manifold of codimension one can be defined and studied 
by Pfaff's nonintegrable (complete integrable) equation $\eta = 0$.

{\em The present paper defines and studies a structure for some $(4n+3)$-dimensi-onal 
manifolds which is named almost coquaternion structure. This structure is composed 
of three almost cocomplex structures $(\phi_a, \xi_a, \eta_a)$, $a = 1,2,3$, 
which satisfy some relations and may be considered as analogous to the 
almost quaternion structure for $(4n+4)$-dimensional manifolds $[49]- [53]$}.

Independently, Y. Kuo [23] introduced in the same way the same structure 
calling it an {\it almost contact 3-structure}. Then, S. Sasaki [36], S. Tachibana and W. Yu [39], 
S. Tanno [40] etc investigated some properties of this structure using the denomination given by Kuo.

The sphere $S^{4n+3}$ is a typical example of differentiable manifold which 
admits an almost coquaternion structure $(\phi_a, \xi_a, \eta_a)$, $a = 1,2,3$ [49], [36]. 
Using the 1-forms $\eta_a$ of the almost coquaternion structure of the sphere $S^{4n+3}$, 
C. Teleman [44] defined and studied on $S^{4n+3}$ a nonholonomic manifold $V^{4n}_{4n+3}$ 
whose Riemannian metric is the one of a symmetric space of E. Cartan [4]. Keeping in mind 
Teleman's idea, we observed that on an almost coquaternion manifold a nonholonomic 
(holonomic) manifold of codimension three can be defined and studied by 
nonintegrable (complete integrable) Pfaff's system $\eta_1 = 0$, $\eta_2 = 0$, $\eta_3 = 0$.

Our results are included in seven Chapters whose summary is given further on.

\newpage

\centerline{\bf CONTENTS}
\vspace{0.5cm}
\centerline{\bf Chapter I: Almost coquaternion structures}
\vspace{0.3cm}
\noindent$\S 1$ Algebraic preliminaries\\
$\S 2$ Almost coquaternion structures\\
$\S 3$ The structural group of the tangent bundle I\\
$\S 4$ Associated Riemannian metric\\
$\S 5$ Almost coquaternion metric structures\\
$\S 6$ The structural group of the tangent bundle II

\vspace{1cm}
\centerline{\bf Chapter II: Almost coquaternion structures and} 
\vspace{-0.05cm}
\centerline{\bf theirs connection with almost quaternion structures}
\vspace{0.3cm}
\noindent$\S 1$ On the geometry of almost quaternion manifolds\\
$\S 2$ The almost coquaternion structure of the manifold $\overline M\times R^3$\\
$\S 3$ The almost coquaternion structure of the manifold $M\times R^3$\\
$\S 4$ Almost quaternion contact manifolds and\\ their product with the almost coquaternion manifolds

\vspace{1cm}
\centerline{\bf Chapter III: Almost coquaternion affine connections} 
\vspace{0.3cm}
\noindent$\S 1$ The bundle of coquaternion linear frames\\
$\S 2$ Almost coquaternion affine connections\\
$\S 3$ $(\xi_a,\eta_a)$-connections\\
$\S 4$ Almost coquaternion metric connections\\
$\S 5$ The holonomy group and the almost coquaternion structure

\vspace{1cm}
\centerline{\bf Chapter IV: Almost coquaternion hypersurfaces} 
\vspace{0.3cm}
\noindent$\S 1$ Almost coquaternion hypersurfaces\\
$\S 2$ Induced $(\phi_a,\xi_a,\eta_a)$-connection\\
$\S 3$ Almost coquaternion metric hypersurface\\
$\S 4$ Imbedding of almost coquaternion Riemannian manifolds\\
$\S 5$ $A_m^0$-manifolds and $\eta$-Einstein manifolds

\vspace{1cm}
\centerline{\bf Chapter V: Fibering of almost coquaternion manifolds} 
\vspace{0.3cm}
\noindent$\S 1$ $Q_xM$-invariant almost coquaternion structures\\
$\S 2$ Fiberings of almost coquaternion manifolds\\
$\S 3$ Fiberings of almost coquaternion Riemannian manifolds\\
$\S 4$The curvature form

\vspace{1cm}
\centerline{\bf Chapter VI: On $\phi$-transformations} 
\vspace{0.3cm}
\noindent$\S 1$ $\phi$-transformations\\
$\S 2$ Automorphisms\\
$\S 3$ The almost coquaternion structure of\\ the $(4n+3)$-dimensional Lie groups

\vspace{1cm}
\centerline{\bf Chapter VII: Examples} 
\vspace{0.3cm}
\noindent$\S 1$ The almost coquaternion structure of $3$-dimensional orientable manifolds\\
$\S 2$ The case of the sphere $S^3$\\
$\S 3$ The case of the Tzitzeica hypersurface $xyzt=1$\\
$\S 4$The sphere $S^{4n+3}$ as an almost coquaternion manifold\\
$\S 5$ The almost coquaternion structure of $T^1\overline M$
\vspace{1cm}

\centerline {\bf References}

\newpage

\centerline{\bf \large CHAPTER I}

\bigskip

We suppose that all the used differentiable manifolds and maps are of class $C^\ty$ and 
we denote by $\mathcal{X}(M)$ the Lie algebra of all vector fields on the manifold $M$.

Let $M$ be a $(4n+3)$-dimensional manifold.

{\bf 1.1. Definition} [49], [23]. {\em An almost coquaternion structure on $M$ is an aggregate 
consisting of three almost cocomplex structures $(\phi_a, \xi_a, \eta_a)$, $a = 1,2,3$, which satisfy
$$
\phi_a \circ \phi_b - \xi_a \otimes \eta_b = - \phi_b \cdot \phi_a + \xi_b \otimes \eta_a = \phi_c
$$
$$
\phi_a \xi_b = - \phi_b \xi_a = \xi_c
$$
$$
\eta_a \cdot \phi_b = - \eta_b \cdot \phi_a = \eta_c,
$$
$$
\eta_a (\xi_b) = \eta_b (\xi_a) = 0,
$$
for any cyclic permutation $\{a,b,c\}$ of $\{1,2,3\}$.}

In the paper [39] it is shown that the aggregate $(\phi_a, \xi_a, \eta_a)$, $a = 1,2,3$, 
cannot be completed with a fourth almost cocomplex structure which, together with 
the three structures given above satisfy relations of the same type as those of 
Definition 1.1. On the other hand, we observe [49], [23] that in order to have an 
almost coquaternion structure on $M$ it is sufficient to give only two almost 
cocomplex structures which satisfy certain relations, since the third structure 
results from the given structures. The more explicit form is [49]:

{\bf 1.2. Theorem}. {\em If a differentiable manifold $M$ admits two almost cocomplex 
structures $(\phi_a, \xi_a, \eta_a)$, $a = 1,2$, which satisfy
$$
\phi_1 \circ \phi_2 + \phi_2 \circ \phi_1 = \xi_1 \otimes \eta_2 + \xi_2 \otimes \eta_1,
$$
$$
\eta_1 (\xi_2) = 0,
$$
then we have $\eta_2 (\xi_1) = 0$, $\phi_1 (\xi_2) = -\phi_2\xi_1$, $\eta_1 \circ \phi_2 = -\eta_2 \circ \phi_1$, 
and $M$ admits a third almost cocomplex structure $(\phi_3, \xi_3, \eta_3)$ defined by
$$
\phi_3 = \phi_1 \circ \phi_2 - \xi_1 \otimes \eta_2 = - \phi_2 \circ \phi_1 + \xi_2 \otimes \eta_1,
$$
$$
\xi_3 = \phi_1 \xi_2 = - \phi_2 \xi_1, \quad  \eta_3 = \eta_1 \circ \phi_2 = - \eta_2 \circ \phi_1,
$$
so that the relations of Definition 1.1 are valid}.

{\bf 1.3. Theorem}. {\em If $(\phi_a, \xi_a, \eta_a)$, $a = 1,2,3$, is an 
almost coquaternion structure on $M$, then $\fo (A^a_d) \in SO (3)$, $\fo\al: M \to R \setminus \{0\}$,
$$
\left(\bar\phi_d = \sum_a A^a_d \phi_a, \; \bar\xi_d = \frac{1}{\al} \sum_a A^a_d \xi_a, \; \bar\eta_d = \al \sum_a A^a_d \eta_a\right), \quad d = 1,2,3,
$$
is again an almost coquaternion structure on $M$}.

$M$ is said to be an {\em almost coquaternion manifold}.

Let $x$ be an arbitrary point of $M$. We denote the tangent space at $x$ by $T_x M$ 
and its hypercomplexification $H \otimes_R T_xM$ by $T^H_x M$, where $H$ is the corps 
of quaternion numbers and $R$ is the real number field.

{\bf 1.4. Theorem}. {\em If by $D_{ax} M$, $a=1,2,3$, we understand unidimensional 
distributions defined respectively by vector fields $\xi_a$, then the existence of an 
almost coquaternion structure on $M$ is equivalent to the existence of the 
disjoint distributions $S^H_x M$, $D_{ax} M$ such that}
$$
T^H_x M = S^H_x M \oplus \tau' S^H_x M \oplus\tau''S^H_x M \oplus \tau'''S^H_x M \oplus D^H_{1x} M \oplus D^H_{2x} M \oplus D^H_{3x} M.
$$

Thus, the existence of an almost coquaternion structure on $M$ implies $\hbox{dim} M = 3n+3$.

Every almost coquaternion manifold $M$ admits an {\em associated Riemannian metric} 
to the given almost coquaternion structure in the sense of

{\bf 1.5. Theorem.} {\em If $M$ is an almost coquaternion manifold, 
then there exists a positive definite Riemannian metric $g$ such that
$$
\eta_a (X) = g(\xi_a, X),
$$
$$
g(\phi_a X, \phi_a Y) = g(X,Y) - \eta_a (X) \eta_a (Y), \; a = 1,2,3, \; \fo X,Y \in \mathcal{X} (M).
$$

The aggregate $(\phi_a, \xi_a, \eta_a, g)$, $a = 1,2,3,$ is called an 
almost coquaternion metric structure and $M$ is said to be an almost coquaternion Riemannian manifold}.

An almost coquaternion metric structure can be described by means of 1-forms $\eta_a$ and 
2-forms $\Theta_a (X,Y) = g(\phi_a X,Y)$. Furthermore, we have

{\bf 1.6. Theorem}. {\em If $(\phi_a, \xi_a, \eta_a, g)$, $a = 1,2,3,$ is an 
almost coquaternion metric structure, then, $\fo (A^a_d) \in SO(3)$, $\fo \al:M \to (0,\ty)$,
$$
\left(\sum_a A^a_d \phi_a, \frac{1}{\al} \sum_a A^a_d \xi_a, \; \al \sum_a A^a_d \eta_a, \; \al g + (\al^2 - \al) \sum_a \eta_a \otimes \eta_a\right), \; d = 1,2,3,
$$
is again an almost coquaternion metric structure on $M$}.

Also, we find

{\bf 1.7. Theorem}. {\em If $M$ is am almost coquaternion  (Riemannian) manifold, 
then the structural group of its tangent bundle $TM$ reduces to $GL(n,H) \ti 1 \ti 1 \ti 1$ $(Sp (n) \ti 1\ti 1 \ti 1)$. 
The converse is also true}.

Consequently, every almost coquaternion manifold $M$ is orientable.

Some of the results presented in this Chapter have been published in [49].
\bigskip

\centerline{\bf \large CHAPTER II}

\bigskip

Let $\overline M$ be a $4n$-dimensional differentiable manifold which possesses an 
almost quaternion Hermitian structure $(J_a, G)$, $a = 1,2,3,$ and $M$ be a 
manifold which admits an almost coquaternion metric structure $(\phi_a, \xi_a, \eta_a, g)$, $a = 1,2,3$.

{\bf 2.1. Theorem}. {\em The product manifold $\overline M \ti M$ possesses an 
almost coquaternion metric structure $(\phi'_a, \xi'_a, \eta'_a, g)$ defined by}
$$
\phi'_a \left(\begin{array}{c} X \\ Y \end{array} \right) =
\left(\begin{array}{c} J_a X \\ \phi_a Y \end{array} \right), \; \fo X \in T_x \overline M, \; \fo Y \in T_y M, \; (x,y) \in \overline M \ti M,
$$
$$
\xi'_a = \left(\begin{array}{cc} 0 \\ \xi_a \end{array} \right), \; \eta'_a = (0,eta_a), \; g' = \left(\begin{array}{cc} G & 0 \\ 0 & g \end{array} \right).
$$

Taking into account that the aggregate $(0,\frac{d}{dt}, dt, 1)$ can 
be regarded as an almost cocomplex metric structure on the real line $R$, we have

{\bf 2.2. Theorem.} {\em If $M$ is an almost coquaternion Riemannian manifold, 
then the product manifold $M \ti R$ has an almost quaternion Hermitian structure
$$
J_a = \left(\begin{array}{cc} \phi_a & \xi_a \\ -\eta_a & 0 \end{array} \right), \quad G = \left(\begin{array}{cc} g & 0 \\ 0 & 1 \end{array} \right)
$$
induced by $(\phi_a, \xi_a, \eta_a, g)$ and $(0,\frac{d}{dt}, dt, 1)$}.

As in [35], we obtain

{\bf 2.3. Theorem}. {\em The almost coquaternion structure determines the following tensor fields on $M$
$$
N^1_a (X,Y) = N_a (X,Y,\phi_a) + 2 d\eta_a (X,Y) \xi_a, \; \hbox{of type} \; (1,2),
$$
$$
N^2_a (X,Y) = -2\eta_a (\phi_a X,Y) - 2 d\eta_a (X,\phi_a Y) \xi_a, \; \hbox{of type} \; (0,2),
$$
$$
N^3_a (X) = -(L_{\xi_a} \phi_a) (X), \; \hbox{of type} \; (1,1),
$$
$$
N^4_a (X) = -(L_{\xi_a} \eta_a) (X), \; \hbox{of type} \; (0,1), \; \fo X,Y \in \mathcal{X} (M),
$$
where $N_a (X,Y; \phi_a)$ is the Nijenhuis tensor of the endomorphism $\phi_a$}.

The tensor fields satisfy some relations which corresponding to those satisfied by the 
Nijenhuis tensors attached respectively to the three almost complex structures 
$J_a$ on $M \ti R$. As an application of these relations, we deduce

{\bf 2.4. Theorem}. {\em $N^1_1, N^1_2, N^1_3$ vanish if any two of them vanish.}

On the other hand, using the structure tensor $T$ of the almost quaternion structure on 
$M \ti R$ we deduce that the almost coquaternion structure determines the following tensor fields on $M$
$$
T^1 (X,Y) = \frac{2}{3} \sum_a N^1_a (X,Y), \quad \hbox{of type} \; (1,2),
$$
$$
T^2 (X,Y) = \frac{2}{3} \sum_a N^2_a (X,Y), \quad \hbox{of type} \; (0,2),
$$
$$
T^3 (X,Y) = \frac{2}{3} \sum_a N^3_a (X,Y), \quad \hbox{of type} \; (1,1),
$$
$$
T^4 (X,Y) = \frac{2}{3} \sum_a N^4_a (X,Y), \quad \hbox{of type} \; (0,1), \; \fo X,Y \in \mathcal{X} (M).
$$
Notice now, that $T^1$, $T^2$, $T^3$, $T^4$ satisfy some relations which correspond to 
those satisfied by the structure tensor $T$ attached to the almost quaternion structure on $M \ti R$. These imply

{\bf 2.5. Theorem.} {\em If any one of $T^1$, $T^2$, $T^3$, vanishes, then $T^4$ vanishes. If $T^1$ 
vanishes, then all the other tensors $T^2$, $T^3$, $T^4$ vanish}.

{\bf 2.6. Theorem.} {\em $T^1$ vanishes if and only if $N^1_1$, $N^1_2$, $N^1_3$ simultaneously vanish}.

From the former discussions, it appears that the almost coquaternion structure 
$(\phi_a, \xi_a, \eta_a)$, $a = 1,2,3,$ on $M$ whose tensor $T^1$ vanishes (identically) 
corresponds to the pseudo-quaternion structure $J_a$, $a=1,2,3,$ on $M \ti R$. 
Hence it consists of three normal almost cocomplex structures.

{\bf 2.7. Definition}. {\em The almost coquaternion structure whose tensor $T^1$ 
vanishes identically is called a pseudo-coquaternion structure and the manifold 
with such a structure a pseudo-coquaternion on manifold.}

The induced almost quaternion structure on $M \ti R$ may be integrable. So, we admit

{\bf 2.8. Definition}. {\em If the almost quaternion structure $J_a$ on $M \ti R$ is integrable, 
then we say that the almost coquaternion structure $(\phi_a, \xi_a, \eta_a)$ on $M$ 
is normal and that $M$ is a normal almost coquaternion manifold.}

A normal almost coquaternion structure is a pseudo-coquaternion structure. 
The converse is not generally true.

Finally, we observe that if $N$ is a $(4m+1)$-dimensional manifold which possesses an 
almost quaternion contact metric structure $(\va_a, \xi, \eta, h)$, $a = 1,2,3$, [16], then we have

{\bf 2.9. Theorem}. {\em The product $M\ti N$ of an almost coquaternion 
Riemannian manifold $M$ and an almost quaternion contact Riemannian manifold 
$N$ is an almost quaternion Hermitian manifold with the structure $(J_a, G)$, $a=1,2,3,$ defined by}
$$
J_a \left(\begin{array}{c} X \\ Y \end{array} \right) = 
\left(\begin{array}{c} \phi_a X + \eta (Y) \xi_a \\ \va_a Y - \eta_a (X) \xi \end{array} \right), \; \fo X \in T_x M, \; \fo Y \in T_y N,
$$
$$
G = \left(\begin{array}{cc} g & 0 \\ 0 & h \end{array} \right).
$$

The results obtained previously are included in [53].

\bigskip

\centerline{\bf \large CHAPTER III}

\bigskip

First we show that the following theorem is valid.

{\bf 3.1. Theorem}. {\em An almost coquaternion (metric) structure on $M$ is 
equivalent to the existence of a $G$-structure $P_G (M)$, 
where $G$ is the real representation of $GL (n,H) \times 1 \times 1 \times1)$.}

This $G$-structure was called the {\em bundle of coquaternion linear frames} 
({\em the bundle of coquaternion orthonormal frames}).

{\bf 3.2. Definition}. A linear - or affine - connection of $M$ is called an almost 
coquaternion (metric) connection if it is determined by a connection in the bundle of 
coquaternion linear frames (the bundle of coquaternion orthonormal frames).

We have

{\bf 3.3. Theorem.} {\em  A linear connection $\nabla$ of an almost coquaternion 
manifold $M$ with almost coquaternion structure $(\phi_a, \xi_a, \eta_a)$, $a = 1,2,3,$ 
is an almost coquaternion connection if and only if 
$\phi_1, \phi_2, \xi_1$ $(\Longrightarrow \phi_3, \xi_2, \xi_3, \eta_a)$ are parallel with respect to $\nabla$}.

{\bf 3.4. Theorem.} {\em $\nabla$ is an almost coquaternion affine connection 
if and only if the connection $1$-forms $\omega^R_S$ of the coquaternion frame fields
$$
X_S: X_s, \; X_{s'} = \phi_1 X_s, \; X_{s''} = \phi_2 X_s, \; X_{s'''} = \phi_3 X_s, \; X_{4n+a} = \xi_a,
$$
$$
S = s, s', s'', s''', 4n+a; \; s = 1,2, \ld, n,
$$
are given by the matrix
$$
\omega = \left( \begin{array}{c|c}
\begin{array}{cccc} \om^r_1 & \om^r_{s'} & \om^r_{s''} & \om^r_{s'''} \\ \noa - \om^r_{s'} & \om^r_{s} & - \om^r_{s'''} & \om^r_{s''} \\ \noa - \om^r_{s''} & \om^r_{s'''} & \om^r_{s} & - \om^r_{s'} \\ \noa -\om^r_{s'''} & -\om^r_{s''} & \om^r_{s'} & \om^r_{s} \end{array} & 0 \\ \noa \hline 0 & \begin{array}{ccc} 0 & 0 & 0 \\ \noa 0 & 0 & 0 \\ \noa 0 & 0 & 0 \end{array}
\end{array} \right).
$$}

In fact $X_S$ constitute a system of $4n+3$ independent congruences in 
$U \su M$ and $\om^R = ds^R$. Therefore, using the connection 1-forms we 
can build up the components $\ga^R_{TS}$ of the coquaternion affine connection 
$\nabla$ with respect to the congruences $X_S$. Taking into account the definition 
of these components [54], [10-14] expressed in the invariant language [48], we find
$$
\ga^R_{TS} = ds^R (\nabla_{X_T} X_S) = \om^R_S (X_T).
$$
So, $\nabla$ is a coquaternion affine connection if and only if its components 
with respect to the congruences $X_S$ are given by $4n+3$ matrices having the form of $\omega$.

We also prove the following theorems.

{\bf 3.5. Theorem.} {\em Every almost coquaternion manifold $M$ admits an affine connection $\hat\na$ such that}:

1) {\em $\xi_a$ and $\eta_a$ are parallel with respect to $\hat\na$},

2) {\em the torsion $\hat T$ of $\hat\na$ is given by
$$
\hat T (X,Y) = 2\sum_a d\eta (X,Y)\xi_a, \quad \fo X,Y,
$$
where $X$ and $Y$ are vector fields}.

{\bf 3.6. Theorem}. {\em A linear connection $\na$ of an almost coquaternion 
Riemannian manifold $M$ with almost coquaternion metric structure 
$(\phi_a, \xi_a, \eta_a, g)$ is an almost coquaternion metric connection 
if and only if $\phi_1, \phi_2, g (\Rightarrow \phi_3, \xi_a, \eta_a)$ are parallel with respect $\na$}.

{\bf 3.7. Theorem}. {\em $\na$ is an almost coquaternion metric connection 
if and only if the connection 1-forms $\om^R_S$ of the coquaternion 
orthonormal frame fields $X_S$ are given by the matrix $\om$ with the additional conditions $\om^R_S = - \om^S_R$}.

Some results concerning holonomy and almost coquaternion (metric) structure on $M$ are given at the end of this chapter.

These results are included in [51].

\bigskip

\centerline{\bf \large CHAPTER IV}

\bigskip

Suppose $\overline M$ $(dim \overline M = 4n+4)$ has an almost quaternion structure 
$J_a$, $a = 1,2,3$, and $M$ is an orientable hypersurface imbedded in 
$\overline M$ by $i: M \to \overline M$. Let $TM$ denote the tangent bundle of 
$M$ and $T_r \overline M$ the restriction of the tangent bundle of $\overline M$ to $M$. 
We denote by $i_*$ the differential of the imbedding $i$ so that $i_*$ is the 
mapping $i_*: TM \to T_R \overline M$. Let $\zeta$ be a pseudonormal vector field 
defined along $M$, i.e., $\zeta \in T_r \overline M$ and $\zeta \notin TM$. Then we can 
find a mapping $i^{-1}_*: T_r \overline M \to TM$ and a 1-form $\om$ defined on $M$ such that
$$
i^{-1}_* \circ i_* = id, \quad i_* \circ i^{-1}_* = Id - \zeta \otimes \om,
$$
$$
\om \circ i_* = 0, \quad i^{-1}_* \zeta = 0, \; \om (\zeta) = 1, \; \om (J_a \zeta) = 0,
$$
where $id.$ is the identity on $TM$ and $Id$. is the identity on $T_r \overline M$.

{\bf 4.1. Theorem}. {\em An orientable hypersurface $M$ of an almost quaternion manifold 
$\overline M$ has a naturally induced almost coquaternion structure
$$
\phi_a = i^{-1}_* \circ J_a \circ i_*, \quad \xi_a = i^{-1}_* J_a \zeta, \quad \eta_a = -\om \circ J_a \circ i_*, \; a = 1,2,3.
$$}

Let $D$ be an affine connection on $\overline M$. The affine connection of 
$\overline M$ induces an affine connection $\na$ and three fundamental 
tensor fields: $h$ of type (0,2) of type (1,1), $m$ of type (0,1), on $M$ in a natural manner [8], [37]. We have

{\bf 4.2. Theorem}. {\em Let $D$ be an almost quaternion connection on $\overline M$, i.e.,
$$
(D_X J_a)(Y) = 0, \quad \fo X,Y \in \mathcal{X} (\overline M).
$$
The connection $\na$ on $M$ is an almost coquaternion connection if and only if $h=0$, $l=0$, $m=0$.}

Suppose $\overline M$ has an Hermitian metric $G$. Then the induced metric $g$ is given by
$$
g(X,Y) = G(i_* X, i_* Y), \quad \fo X,Y \in T_xM.
$$
$\overline M$ is orientable and we assume that $\overline M$ is oriented. Then we can 
choose a differentiable field $\zeta$ of unit normal vectors over $M$ if and only if $M$ is orientable. 
The field $\zeta$ satisfies the equations $G(\zeta, \zeta) = 1$, $G(\zeta, i_* X) = 0$, $\fo X \in T_xM$ 
and $\om  (J_a \zeta) = 0$, $a = 1,2,3,$ together with its associated 1-form $\om (X) = G(\zeta, X)$.

{\bf 4.3. Theorem}. {\em An orientable hypersurface $M$ of an almost quaternion Hermitian 
manifold $\overline M$ has a naturally induced almost coquaternion metric structure}
$$
\phi_a = i^{-1}_* \circ J_a \circ i_*, \quad \xi_a = i^{-1}_* J_a \zeta, \quad \eta_a = -\om \circ J_a \circ i_*, \; a = 1,2,3.
$$
$$
g(X,Y) = G(i_* X, i_* Y), \; \fo X,Y \in T_xM.
$$

Particularly, any hypersurfaces in $R^{4n+4}$ which has one of the properties:

1) is simply connected and connected,

2) can be defined implicitly by $M: f(x) = c$, where the differential $df$ is not zero at any point of $M$,

3) is convex,
\newline
is orientable and hence it possesses an almost coquaternion metric structure. 
For example, the sphere $S^{4n+3}$ is a compact manifold which has an almost coquaternion 
metric structure and the hypersurface of Tzitzeica $x^1x^2 \ld x^{4n+4} = 1$ is a 
noncompact manifold which has an almost coquaternion metric structure.

If $R$ is the real line, then we have

{\bf 4.4. Theorem}. {\em An almost coquaternion Riemannian manifold $M$ can be 
imbedded in the almost quaternion Hermitian manifold $M \ti R$ as a totally umbilical or geodesic hypersurface.}

The above mentioned ideas arose from [41] and were published in [50].

Let $M$ be a $m$-dimensional manifold and $\Gamma^i_{jk}$, $i,j,k,l = 1, 2, \ld, m$, 
be a connection on $M$. We denote the Ricci tensor field attached to the given connection by $\Gamma_{ij}$.

{\bf 4.5. Definition} [15]. {\em If $\na_k (\det (\Gamma_{ij})) = 0$, then the manifold $M$ is called an $A^0_m$-manifold.}

{\bf 4.6. Theorem}. {\em If $\det (\Gamma_{ij}) \ne 0$, then $M$ is an $A^0_m$-manifold if and only if
$$
\Gamma^{ij} \na_k \Gamma_{ij} - \Gamma^{ij} \Gamma^{l}_{jk} \Gamma_{ei} + \Gamma^l_{lk} = 0, 
\quad \Gamma^{ij} \stackrel{def}{=} \frac{1}{\det (\Gamma_{ij})} \frac{\pa \det (\Gamma_{ij})}{\pa \Gamma_{ij}},
$$
or
$$
\Gamma^{ij} \na_k \Gamma_{ij} = 0,
$$
if $\Gamma_{ij}$ is a symmetrical tensor}.

Using this Theorem we have:

{\bf 4.7. Theorem}. {\em Any $\eta$-Einstein manifold is an $A^0_m$-manifold.}

\bigskip

\centerline{\bf \large CHAPTER V}

\bigskip

Being given an almost coquaternion manifold $M$ with structure 
$(\phi_a, \xi_a, \eta_a)$, $a = 1,2,3,$ we suppose that $\{\xi_a\}$ gives rise 
to a Lie group of transformations on $M$. Let us denote by $G$ this 3-dimensional group.

{\bf 5.1. Definition}. {\em An almost coquaternion structure $(\phi_a, \xi_a, \eta_a)$ is called $Q_x M$-invariant if the distribution
$$
Q_xM = \{X| X \in T_xM, \; \eta_1 (X) = 0, \eta_2 (X) = 0, \eta_3 (X) = 0\}
$$
is invariant with respect to $R_t$, where $R_t$ means the right translation in $M$ by $t\in G$.}

{\bf 5.2. Theorem}. {\em An almost coquaternion structure is $Q_xM$-invariant if and only if
$$
\begin{array}{l}
L_{\xi_1} \eta_a = C^a_{21} \eta_2 + C^a_{31} \eta_3, \\ \noa
L_{\xi_2} \eta_a = C^a_{12} \eta_1 + C^a_{32} \eta_3, \\ \noa
L_{\xi_3} \eta_a = C^a_{13} \eta_1 + C^a_{23} \eta_2, \end{array}
$$
where $C^a_{bc}$ are the structure constants of the Lie algebra 
of $G$ and $L_{\xi_a}$ represents the Lie differentiation with respect to $\xi_a$}.

{\bf 5.3. Lemma}. {\em If the structure $(\phi_a, \xi_a, \eta_a)$ is $Q_xM$-invariant, 
then the tensor field $\Gamma = \di\sum_a \xi_a \otimes \eta_a$ is invariant under the action of $G$. The converse is also true}.

As the tensor field of type (2,2)
$$
P = \frac{1}{4} (id \otimes id + \sum_a \phi_a \otimes \phi_a + \Gamma \otimes\Gamma)
$$
is invariant by $G$ if and only if the field $\di\sum_a \phi_a \otimes \phi_a$ invariant by $G$, we admit

{\bf 5.4. Definition}. {\em A $Q_xM$-invariant almost coquaternion structure is said to 
be $P$-invariant if and only if $\di\sum_a \phi_a \otimes \phi_a$ is invariant under the action of $G$.}

Let us suppose that $M$ is an almost coquaternion Riemannian manifold with structure 
$(\phi_a, \xi_a, \eta_a,g)$, $a = 1,2,3,$ and that $\{\xi_a\}$ determine a Lie group of 
motions with respect to $g$. In these conditions we get

{\bf 5.5. Theorem}. {\em The group $G$ is necessarily isomorphic to a unitary, 
semi-simple group, if $C^1_{23} = m\ne 0$, or is isomorphic to an Abelian group, 
if $C^1_{23} = 0$. If $G$ is isomorphic to a unitary, semi-simple group, then its 
representative Riemanian space has a positive constant curvature}.

We specify that in this case the structure $(\phi_a, \xi_a, \eta_a)$ is $Q_xM$-invariant as
$$
L_{\xi_a} \eta_a = 0, \qu L_{\xi_a} \eta_b = - L_{\xi_b} \eta_a = m \eta_c,
$$
for any cyclic permutation $\{a,b,c\}$ of $\{1,2,3\}$. Evidently, the tensor field 
$\Gamma = \di\sum_a \xi_a \otimes \eta_a$ is invariant with respect to such a group $G$ and we have

{\bf 5.6. Lemma}. {\em If $\{\xi_a\}$ generates a unitary, semi-simple group of motions, with respect to $g$, such that
$$
L_{\xi_a} \phi_a = 0, \qu L_{\xi_a} \phi_b = - L_{\xi_b} \phi_a = m \phi_c,
$$
for any cyclic permutation $\{a,b,c\}$ of $\{1,2,3\}$, then the tensor fields
$$
\begin{array}{ll}
\sum_a \phi_a \otimes \phi_a, & \hbox{of type} \; (2,2), \\ \noa
P = \di\frac{1}{4} \left(id \otimes id + \di\sum_a \phi_a \otimes \phi_a + \Gamma \otimes\Gamma\right), & \hbox{of type} \; (2,2), \\ \noa \phi_1 \otimes \phi_2 \otimes \phi_3, & \hbox{of type} \; (3,3), \\ \noa \di\sum_a \xi_a \otimes \phi_a, & \hbox{of type} \; (2,1), \end{array}
$$
are invariant under the action of $G$}.

Let $M$ be an almost coquaternion manifold with structure $(\phi_a, \xi_a, \eta_a)$, $a = 1,2,3$. We admit

{\bf 5.7. Definition}. {\em If $\{\xi_a\}$ gives rise to a Lie group of transformations on $M$ and the distribution defined by $\{\xi_a\}$ id regular in Palais' sense [32], then the almost coquaternion structure $(\phi_a, \xi_a, \eta_a)$ is said to be strictly regular.}

Suppose $M$ is a connected compact manifold. We denote by $M_G$ the 
quotient space and by $\pi: M \to M_G$ the canonical projection. It is known that 
[32] $M_G$ is a connected compact manifold for which $\pi (U) (x^1, \ld, x^{4n})$ is a coordinate neighborhood.

{\bf 5.8. Theorem}. {\em If $(\phi_a, \xi_a, \eta_a)$ is a strictly regular and 
$Q_xM$-invariant almost coquaternion structure on a connected compact manifold $M$, then}

(i) {\em $M$ is a principal $G$-bundle over $M_G$ whose fiber is compact} [32],

(ii) {\em $\Gamma = \di\sum_a \xi_a \otimes \eta_a$ is a connection on $M$, and}

(iii) {\em $\eta = \{\eta_a\}$ is a connection 1-form on $M$}.

The Theorem 5.8 shows that we may consider principal bundles of type $M(M_G, G, \pi)$, 
where the total space $M$ of possesses a $Q_xM$-invariant almost coquaternion 
structure and $G$ is the group generated by $\{\xi_a\}$.

{\bf 5.9. Theorem}. {\em If $M(M_G,G,\pi)$ is a principal bundle in which the total 
space $M$ possesses a $P$-invariant structure $(\phi_a, \xi_a, \eta_a)$ then the tensor field
$$
P = \di\frac{1}{4} \left(id \otimes id + \di\sum_a \phi_a \otimes \phi_a + \Gamma \otimes\Gamma\right)
$$
induces on $M_G$ a global projector on the tensors of degree 2,
$$
\mathcal{P}_p (X,Y) = d\pi (P_x)(X^*_x, Y^*_x),
$$
$\fo p \in M_G$, $x\in M$, $p = \pi (x)$, $X^*_x$ denotes the lift of $X \in \mathcal{X} (M_G)$ with respect to the connection $\eta$}.

Suppose now that $M$ possesses an almost coquaternion metric structure 
$(\phi_a, \xi_a, \eta_a,g)$, $a = 1,2,3$. Extending a result from [49] we get

{\bf 5.10. Theorem}. {\em If the connected compact manifold $M$ has an 
almost coquaternion metric structure for which}: (1) $\Theta_a = d\eta_a$, 
(2) {\em the induced almost quaternion structure on $M \ti R$ is a pseudo-quaternion structure}, 
(3) {\em the distribution defined by $\{\xi_a\}$ is regular in Palais' sense} [32], {\em then}

(i) $M$ {\em is a principal bundle whose fiber is a sphere $S^3$ or a real projective space $RP^3 = S^3 / \{id., -id\}$.}

(ii) $\Gamma = \di\sum_a \xi_a \otimes \eta_a$ {\em is a connection on $M$, and}

(iii) {\em $\eta = \{\eta_a\}$ is a connection 1-form on $M$}.

As $g$ is a $G$-invariant Riemannian metric on $M$, we can define a Riemannian metric on $M_g$ by
$$
h(X,Y) = g(X^*, Y^*), \quad X,Y \in \mathcal{X} (M_G),
$$
called the {\it induced Riemannian metric}.

Next we suppose that the hypotheses of Theorem 5.10 are satisfied and hence Lemma 5.6 is valid. 
Therefore all the tensor fields mentioned in lemma 5.6 induce tensor fields on $M_G$ whose properties are included in some theorems.

The Theorems stated in this part of the paper hold good for any principal bundle $M(M_G, G, \pi)$ 
in which the total space $M$ possesses an almost coquaternion metric structure and $D$ is a 
unitary semi-simple group of motion satisfying Lemma 5.6.

Finally, we make some observations on the structure equation of the fibering in Theorem 5.10.

\bigskip

\centerline{\bf \large CHAPTER VI}

\bigskip

Let $M$ be an almost cocomplex manifold with structure $(\phi, \xi,\eta)$ and $f$ a transformation of $M$.

{\bf 6.1. Definition}. {\em A transformation $f$ of $M$ which leaves $\phi$ invariant, i.e.,
$$
f_* \circ \phi = \phi \circ f_*
$$
is called a $\phi$-transformation [35].}

{\bf 6.2. Lemma}. {\em If $f$ is a $\phi$-transformation of $M$, then there exists a function $\ep: M \to R$ such that
$$
f_* \xi = \ep \xi, \; f^* \eta = \ep \circ f\eta.
$$}

Consider a subset of infinitesimal $\phi$-transformations $X$ such that there is a system of differential equations
$$
\eta_k \xi^i \frac{\pa^2 X^k}{\pa x^i \pa x^j} = h^i_{kj} (x) \frac{\pa X^k}{\pa x^i} + h_{kj} (x) X^k + h_j (x) \leqno (*)
$$
defined in a neighborhood  of $x \in M$ and satisfied by components $X^i$ of $X$; 
denote by $\mathcal{T}_0 (M;\phi)$ the corresponding subset of $\phi$-transformations 
in the set of all $\phi$-transformations over $M$.

Sasaki's papers [35] suggested us to use some results of [1], [32] in order to prove the following

{\bf 6.3. Theorem}. {\em If $M$ is a compact manifold, then $\mathcal{T}_0 (M;\phi)$ is a Lie group of transformations.}

Let $M$ be an almost coquaternion manifold with structure tensors $(\phi_a, \xi_a, \eta_a)$, $a = 1,2,3$.

{\bf 6.4. Definition}. {\em A transformation $f$ on $M$ which leaves $\phi_1$ and $\phi_2$ invariant, i.e.,
$$
f_* \circ \phi_a = \phi_a \circ f_*, \; a = 1,2,
$$
is called a $(\phi_1, \phi_2)$-transformation.}

{\bf 6.5. Lemma}. {\em If $f$ is a $(\phi_1, \phi_2)$-transformation of $M$, then
$$
f_* \circ \phi_3 = \phi_3 \circ f_*
$$
and there exists a function $\ep: M \to R$ such that}
$$
f_* \xi_a = \ep \xi_a, \; f^* \eta_a = \ep \circ f \eta_a, \; a = 1,2,3.
$$

Let $\mathcal{T}_0 (M;\phi_1, \phi_2)$ be the subset of all $(\phi_1, \phi_2)$-transformations 
which verify the relation (*). We have

{\bf 6.6. Corollary}. {\em If $M$ is a compact almost coquaternion manifold, 
then $\mathcal{T}_0 (M;\phi_1, \phi_2)$ is a Lie group of transformations}.

Let $M$ be an almost cocomplex manifold.

{\bf 6.7. Definition}. {\em A transformation $f$ of $M$ which leaves the structure $(\phi, \xi, \eta)$ invariant, i.e.,
$$
f_* \circ \phi = \phi \circ f_*, \; f_* \xi = \xi \; (\Rightarrow f^* \eta = \eta)
$$
is called an automorphism of $M$ [35].}

It results: if $M$ is a compact manifold, then the set of all automorphisms of $M$ is a Lie group [35].

Let $M$ be an almost coquaternion manifold.

{\bf 6.8. Definition}.  {\em A transformation $f$ of $M$ which leaves the structure $(\phi_a, \xi_a, \eta_a)$ invariant, i.e.,
$$
f_* \circ \phi_1 = \phi_1 \circ f_*, \; f_* \circ \phi_2 = \phi_2 \circ f_*, \; f_* \xi_1 = \xi_1,
$$
$$
(\Rightarrow f_* \circ \phi_3 = \phi_3 \circ f_*, \; f_* \xi_2 = \xi_2, \; f_* \xi_3 = \xi_3, \; f^* \eta_a = \eta_a , \; a = 1,2,3),
$$
is called an automorphism of $M$.}

It follows that the set of all automorphisms of a compact almost coquaternion manifold $M$ is a Lie group.

Finally, we prove

{\bf 6.9. Theorem}. {\em Every $(4n+3)$-dimensional Lie group admits a left invariant almost coquaternion structure.}

The results obtained in this chapter have been included in [52].

\bigskip

\centerline{\bf \large CHAPTER VII}

\bigskip

Let $M$ be a 3-dimensional manifold. We have

{\bf 7.1. Theorem}. {\em $M$ has an almost coquaternion (metric) structure if and only if it is parallelizable.}

The hypothesis that $M$ is parallelizable is equivalent to the fact that it 
possesses three vector fields $\xi_a$ which are linearly independent at 
every point of $M$. Let $\eta_a$ be the dual 1-forms. Define three linear endomorphisms $\phi_a$ by
$$
\phi_1 = \xi_3 \otimes \eta_2 - \xi_2 \otimes \eta_3, \; \phi_2 = \xi_1 \otimes \eta_3 - \xi_3 \otimes \eta_1, \; 
\phi_3 = \xi_2 \otimes \eta_1 - \xi_1 \otimes \eta_2
$$
and a metric tensor $g$ by
$$
g = \eta_1 \otimes \eta_1 + \eta_2 \otimes \eta_2 + \eta_3 \otimes \eta_3.
$$
$(\phi_a, \xi_a, \eta_a, g)$, $a = 1,2,3$, defines an almost coquaternion metric structure on $M$.

As any orientable 3-dimensional manifold is parallelizable [38], we have

{\bf 7.2. Theorem}. {\em Every 3-dimensional orientable manifold has an almost coquaternion metric structure.}

The above presented results permit us to give explicitely almost coquaternion 
metric structures on the sphere $S^3: x^2+y^2+z^2+t^2=1$ and on the Tzitzeica hypersurface $xyzt = 1$.

Then it is shown that the sphere $S^{4n+3}$ possesses a normal 
coquaternion metric structure $(\phi_a, \xi_a, \eta_a, g)$, $a = 1,2,3$, 
for which $\Theta_a = d\eta_a$ and the distribution defined by $\{\xi_a\}$ 
is regular in Palais' sense. So the sphere $S^{4n+3}$ is an example of our 
results obtained in Chapter V. The quotient space of this principal bundle 
is the quaternionic projective space [2], [54].

Finally, we construct some almost coquaternion structures on some 
hypersurfaces of the tangent bundle of an almost quaternion (Hermitian) manifold.

\bigskip

The author expresses his heartiest thanks to Professor Gheorghe Th. Gheorghiu 
for his kind remarks and valuable advices in preparing the PhD Thesis.


\begin{thebibliography}{59}

\bibitem{[1]} Bochner, S., {\em Tensor fields with finite bases}, Ann. of Math. (2), 53 (1951), pp. 400-411.

\bibitem{[2]} Bonan, E., {\em $G$-structures de type quaternionien}, Th\`ese, Paris (1967).

\bibitem{[3]} Bouzon, J., {\em Structures presque-cocomplexes}, Th\`ese, Paris (1965).

\bibitem{[4]} Cartan, E., {\em Sur certaines formes Riemanniennes remarquables des 
g\'eometrie a group fondamental simple}, Annalles de l'Ecole Normale Sup., 44 (1927), pp. 348.

\bibitem{[5]} Chern, S., {\em The geometry of $G$-structures}, Bull. Amer. Math. Soc., 72 (1966), pp. 167-219.

\bibitem{[6]} Clark, R.S. and Bruckheimer, M., {\em Tensor structures on a differentiable manifold}, 
Ann. Mat. Pura Appl. (4), 54 (1961), pp. 123-141.

\bibitem{[7]} Ehresmann, C., {\em Sur la th\'eorie des espaces fibr\'es}, Colloque de Topologie Algebrique, C.N.R.S., Paris (1947), pp. 3.

\bibitem{[8]} Eisenhart, L.P., {\em Non-Riemanian Geometry}, Amer. Math. Soc., New York (1958).

\bibitem{[9]} Gheorghiu, Th. Gh., {\em Hipersuprafe\c te \c Ti\c teica}, Luc\u arile \c stiin\c tifice ale 
Institutului Pedagogic Timi\c soara (1959).

\bibitem{[10]} Gheorghiu, Th. Gh., {\em Sur un calcul tensoriel mixte}, Atti della Acc. Sci. di Bologna, 
Rendiconti, Serie XII, t. III (1965-1966).

\bibitem{[11]} Gheorghiu, Th. Gh., {\em Sur les courbes g\'eod\'esique d'un calcul absolute mixte}, 
Analele St. Univ. Ia\c si, Matematic\u a, t. XI-B (1965), pp. 273-281.

\bibitem{[12]} Gheorghiu, Th. Gh., {\em Connexions proj\'ectives dans un calcul absolute mixte}, 
Rend. Sem. Mat. Messina, t. XI (1966-1967), pp. 21-32.

\bibitem{[13]} Gheorghiu, Th. Gh., {\em Des coordon\'ees normales dans un calcul absolut mixte 
relativement a un \'espace Riemannian}, Rev. Roum. Math. Pures Appl. 13 (1968), pp. 451-468.

\bibitem{[14]} Gheorghiu, Th. Gh., {\em Des coordon\'ees normales  mixtes dans un \'espace a connexion affine}, 
Rend. Sem. Mat. Fis. di Milano XXXVIII (1968).

\bibitem{[15]} Gheorghiu, Th. Gh., {\em O clas\u a particular\u a de spa\c tii cu conexiune afin\u a}, 
St. Cerc. Mat. t. 21, 8 (1969), pp. 1157-1168.

\bibitem{[16]} Hashimoto, S., {\em On differentiable manifolds with almost quaternion contact structure}, 
Tensor (N.S.) 15 (1964), pp. 258-268.

\bibitem{[17]} Ianu\c s, S., \c si Udri\c ste, C., {\em Asupra spa\c tiului fibrat tangent al 
unei variet\u a\c ti diferen\c tiabile}, St. Cerc. Mat. t. 22, 4 (1970), pp. 599-611.

\bibitem{[18]} Ianu\c s, S., {\em Quelques remarques sur l'espace fibr\'e tangent 
d'une vari\'et\'e cosymplectique}, An. St. Univ. Ia\c si, Matematic\u a (1971).

\bibitem{[19]} Ianu\c s, S., {\em Some almost product structure on manifolds with linear connections}, Kodai Math. Sem. Rep. (1971).

\bibitem{[20]} Kobayashi, S. and Nomizu, K., {\em Foundations of differential geometry}, 
vol. I, II, Interscience, New York (1963) and (1969).

\bibitem{[21]} Koszul, J.L., {\em Lectures on fibre bundles and differential geometry}, 
Tata Inst. Fund. Research, Bombay (1960).

\bibitem{[22]} Kraines, V.Y., {\em Topology of quaternionic manifolds}, 
Bull. Amer. Math. Soc. 71 (1965), pp. 526-527; same title, Trans. Amer. Math. Soc. 122 (1966), pp. 357-367.

\bibitem{[23]} Kuo, Y.Y., {\em On almost contact 3-structure}, T\^ohoku Math J. 22 (1970), pp. 325-332.

\bibitem{[24]} Libermann, P., {\em Sur les structures presque complexes et autres 
structures infinit\'esimales r\'eguli\`eres}, Bull. Soc. Math. France 83 (1955).

\bibitem{[25]} Lichnerowicz, A., {\em Th\'eorie globale des connexions et des groupes d'holonomie}, Rome (1955).

\bibitem{[26]} Martinelli, E., {\em Varieta a struttura quaternionale generalizzata}, 
Atti Accad. Naz. Lincei Rend. Cl. Sci. Fis. Mat. Nat. 26 (1959), pp. 353-362.

\bibitem{[27]} Martinelli, E., {\em Modello metrico reale dello spazio proiettivo quaternionale}, 
An. Mat. Pura Appl. 49 (1960), pp. 78-89.

\bibitem{[28]} Obata, M., {\em Affine connections on manifolds with almost complex, quaternion or Hermitian structure}, 
Jap. J. Math. 26 (1956), pp. 43-77.

\bibitem{[29]} Obata, M., {\em Hermitian manifolds with quaternion structure}, T\^ohoku Math. J. 10 (1958).

\bibitem{[30]} Ogiue, K., {\em $G$-structure defined by tensor fields}, K\"odai Math. Sem. Rep. 20 (1968), pp. 54-76.

\bibitem{[31]} Ogiue, K. and Okumura, M., {\em On cocomplex structures}, K\"odai Math. Sem. Rep. 19 (1967), pp. 507-512.

\bibitem{[32]} Palais, R.S., {\em A global formulation of the Lie theory of transformations groups}, 
Mem. Amer. Math. Soc. 22 (1957).

\bibitem{[33]} Papuc, I.D. \c si Albu, C.A., {\em Elemente de geometrie diferen\c tial\u a global\u a}, 
Tipografia Univ. Timi\c soara (1970).

\bibitem{[34]} Sasaki, S., {\em On differentiable manifolds with certain structures which are 
closely related to almost contact structures}, T\^ohoku Math. J. 12 (1960), pp. 459-476; II (with Hatakeyama, Y.) 13 (1961), pp. 281-294.

\bibitem{[35]} Sasaki, S., {\em Almost contact manifolds}, Math. Inst. T\^ohoku Univ. Part I(1965), Part II (1967), Part III (1968).

\bibitem{[36]} Sasaki, S., {\em On spherical space forms with normal contact metric 3-structure}, manuscripts.

\bibitem{[37]} Schouten, J.A., {\em Ricci-Calculus}, Berlin (1954).

\bibitem{[38]} Stiefel, E., {\em Richtungsfelder und fernparallelisms in $n$-dimensionalen Mannigfaltigkeiten}, 
Comment. Math. Helv. 8 (1936), pp. 305-353.

\bibitem{[39]} Tachibana, S. and Yu, W.N., {\em On a Riemannian space admitting more than one 
Sasakian structures}, T\^ohoku Math. J. 22 (1970), pp. 535-540.

\bibitem{[40]} Tanno, S., {\em On the isometry groups of Sasakian manifolds}, 
J. Math. Soc. Japan, 22 (1970), pp. 579-590.

\bibitem{[41]} Tashiro, Y., {\em On contact structure of hypersurfaces in complex manifolds I}, 
T\^ohoku Math. J. 15 (1963), pp. 62-78.

\bibitem{[42]} Tashiro, Y., {\em On contact structure of tangent sphere bundles}, 
T\^ohoku Math. J. 21 (1969), pp. 117-143.

\bibitem{[43]} Teleman, C., {\em Asupra grupurilor de mi\c scare maxime ale 
spa\c tiilor riemanniene $V_n$}, St. Cerc. Mat. 5(1954), pp. 143-171.

\bibitem{[44]} Teleman, C., {\em Sur certaines \'espaces sym\'etriques}, 
Bulletin Scientifique, Acad. R.P.R. t.7, 4 (1955), pp. 977-1002.

\bibitem{[45]} Ti\c teica, G., {\em Sur quelques propri\'et\'es affines}, 
Bull. Math. Phys. de l'Ec. Pol., Bucure\c sti, t. 6(1936), pp. 8-10.

\bibitem{[46]} Ti\c teica, G., {\em Opere}, vol. I, Bucure\c sti (1941).

\bibitem{[47]} Udri\c ste, C., {\em Introducerea unor invarian\c ti diferen\c tiali 
intrinseci intr-o varietate riemannian\u a}, St. Cerc. Mat. 21 (1969), pp. 509-514.

\bibitem{[48]} Udri\c ste, C., {\em Congruences on the tangent bundle of a differentiable manifold}, 
Rev. Roum. Math. Pures Appl., t. 15 (1970), pp. 1079-1096.

\bibitem{[49]} Udri\c ste, C., {\em Structures presque coquaterniennes}, 
Bull. Math. de la Soc. Sci. math. de la R.S.R., t. 13 (61), 4(1969), pp. 487-507. 
O parte din aceast\u a lucrare a fost prezentat\u a la "A' IV-lea congres al matematicienilor 
de expresie latin\u a", Bucure\c sti, 17-24 septembrie 1969.

\bibitem{[50]} Udri\c ste, C., {\em Almost coquaternion hypersurfaces}, 
Rev. Roum. Math. Pures Appl., t. 15, 9(1970), pp. 1545-1553.

\bibitem{[51]} Udri\c ste, C., {\em Coquaternion structures and almost coquaternion structures}, 
Rev. Roum. Math. Pures Appl. (sub tipar).

\bibitem{[52]} Udri\c ste, C., {\em On $\phi$-transformations}, An. St. Univ. Ia\c si, Matematic\u a (1971).

\bibitem{[53]} Udri\c ste, C., {\em On almost coquaternion structures}, 
Studia Univ. Babe\c s-Bolyai, series Math. Phys., Cluj (1971).

\bibitem{[54]} Vr\u anceanu, Gh., {\em Lec\c tii de geometrie diferen\c tial\u a}, vol. I, II (1964).

\bibitem{[55]} Vr\u anceanu, Gh., {\em Asupra unei clase de spa\c tii riemanniene omogene}, 
St. Cerc. Mat., t. 5, 1-2 (1954), pp. 173-223.

\bibitem{[56]} Vr\u anceanu, Gh., {\em Varieta differenziabili anolonome}, 
Rend. Sem. Mat. Messina, t. 10 (1965-1966), pp. 139-146.

\bibitem{[57]} Wakakuwo, H., {\em On Riemannian manifolds with homogeneous 
holonomy groups $Sp(n)$}, T\^ohoku Math. J. 10 (1958).

\bibitem{[58]} Yano, K., {\em On a structure defined by a tensor field $f$ of type $(1,1)$ 
satisfying $f^3 + f = 0$}, Tensor, 14 (1965), pp. 99-109.

\bibitem{[59]} Yano, K., {\em Differential geometry on complex and almost complex spaces}, Pergamon Press (1965).

\end{thebibliography}
\end{document}